\documentclass[reqno,11pt]{amsart}
\usepackage{amssymb}
\usepackage{euscript}

\newcommand{\Rad}{\operatorname{Rad}}
\newcommand{\Ad}{\operatorname{Ad}}
\newcommand{\Aut}{\operatorname{Aut}}
\newcommand{\Lie}{\EuScript}
\newcommand{\prodsemi}{\rtimes}
\newcommand{\iso}{\cong}

\newcommand{\origpaper}{[{\em American Journal of Mathematics} {\bf 109} (1987)
927--961]}

\renewcommand{\=}[1]{\mathrel{\mathop{=}\limits^{(\ref{#1})}}}

\theoremstyle{plain}
\newtheorem{correction}{Proposition 6.5\ignorespaces}

\theoremstyle{definition}
\newtheorem{rmk}{Remark\ignorespaces}	

\newcommand{\refjour}{} %% make sure the name isn't already used

\long\def\refjour[#1]#2 #3: #4. #5.. #6 (#7) #8 \par
 {\bibitem[#1]{#2} #3: #4. {\em#5} {\bf #6} (#7) #8\par}

\begin{document}

\title[Zero-entropy affine maps]{Correction to
\\ ``Zero-entropy affine maps 
\\ on homogeneous spaces"}

\author{Dave Witte}

\address{Department of Mathematics, Oklahoma State University, 
 Stillwater, OK 74078}
 \email{dwitte@math.okstate.edu}

\begin{abstract}
 Proposition~6.4 of the author's paper {\origpaper} is incorrect.
This invalid proposition was used in the proof of Corollary~6.5, so
we provide a new proof of the latter result.
 \end{abstract}

\thanks{I would like to thank Alexander Starkov for bringing this error to my
attention, and for the gracious manner in which he did so.
 This research was partially supported by a grant from the National
Science Foundation.}

\maketitle

Professor A.~Starkov has pointed out an error in the proof of Proposition~6.4
of the author's paper {\origpaper}. Contrary to the assertion near the end of
the first paragraph, it may not be possible to choose $T^{*}$ to be a subgroup
of $\langle g\rangle ^{*}$. (The problem is that $\langle g\rangle ^{*}$ may not
contain a maximal torus of $(\Rad G)^{*}$, because the maximal torus of
$\langle g\rangle ^{*}$ may be diagonally embedded in
$\mathord{\hbox{LEVI}}^{*}\times T^{*}$.) The proposition cannot be salvaged,
so the claim must be retracted.

Fortunately, Proposition~6.4 was used only in the proof of Corollary~6.5, for
which we can give a direct proof. As it is no longer a corollary, we now
reclassify this as a proposition.

\begin{correction}
 Suppose $g$ is an ergodic translation on a locally faithful finite-volume
homogeneous space $\Gamma\backslash G$ of a Lie group~$G$, and assume $G =
\Gamma G^{\circ} = G^{\circ}\langle g\rangle $.
 If $g$ has zero entropy, then, for some nonzero power~$g^{n}$ of~$g$, there is
a finite-volume homogenenous space $\Gamma'\backslash G'$ of some Lie
group~$G'$ whose radical is nilpotent, and a continuous map $\psi\colon
\Gamma'\backslash G' \to \Gamma\backslash G$ that is affine for some
translation $g'\in G'$ via~$g^{n}$. Furthermore, every fiber of~$\psi$ is
finite.
 \end{correction}

\begin{proof} Assume for simplicity that $G$ is connected and simply connected.
(A remark on the general case follows the proof.)  Because $g$ is ergodic, we
may assume $\Gamma\langle g\rangle$ is dense in~$G$.

Let $G^{*}$ be the identity component of the Zariski closure of $\Ad G$ in
$\Aut(\Lie G)$, let $S^{*}$ be a maximal compact torus of the Zariski closure
of $\Ad_{G}\langle g\rangle$, and let $L^{*} T^{*}$ be a reductive Levi
subgroup of~$G^{*}$, containing~$S^{*}$. (So $L^{*}$ is semisimple and $T^{*}$
is a maximal torus of $\Rad G^{*}$ that centralizes~$L^{*}$. From Prop.~6.2, we
know that $T^{*}$ is compact.) The composition of $\Ad_{G}$ and the projection
from~$G^{*}$ onto $T^{*}/(L^{*}\cap T^{*})$ is a homomorphism, which, because
$G$ is simply connected, can be lifted to a homomorphism $\pi\colon G \to
T^{*}$.  Define a map $\phi\colon G \to G \prodsemi T^{*}$ by $x^{\phi} =
(x,x^{-\pi})$, where $x^{-\pi} = (x^{-1})^{\pi}$. Notice that
 \begin{equation} \label{phi-crossed}
 (xy)^{\phi} = x^{\phi}\cdot (y^{x^{-\pi}})^{\phi}, 
 \end{equation}
 for all $x,y\in G$.

The definition of~$\phi$ is based on the nilshadow construction of Auslander
and Tolimieri \cite{Aus-Tol} (or see \cite[\S4]{Witte-super}). In particular,
 $\Rad G^{\phi}$ is the nilshadow of $\Rad G$, so $\Rad G^{\phi}$ is nilpotent.

Let $\overline \Gamma\backslash\overline G$ be the faithful version of
$\Gamma^{\phi}\backslash G^{\phi}$. More precisely, let $\overline \Gamma=
\Gamma^{\phi}/N$ and $\overline G = G^{\phi}/N$, where $N$ is the largest
normal subgroup of~$G^{\phi}$ contained in~$\Gamma^{\phi}$. We know
$\Gamma^{\pi}$ is finite (see Prop.~4.20), so, replacing $\Gamma$ by a
finite-index subgroup, we may assume $\Gamma^{\pi} = e$. This implies that 
 $(\gamma x)^{\phi} = \gamma^{\phi} x^{\phi}$ for all $\gamma\in \Gamma$ and $x
\in G$, so $\phi$ induces a well-defined homeomorphism $\overline \phi\colon
\Gamma\backslash G \to \overline \Gamma \backslash \overline G$. 

Unfortunately, $\overline \phi$~is not affine for~$g$ if $T^{*}$~does not
centralize~$g$. We will compensate for the action of~$T^{*}$ by composing with
a twisted affine map.  Assume for simplicity that $L^{*}\cap T^{*} = e$. (Under
this assumption, the map $L^{*}\times T^{*} \to L^{*}T^{*}$ is an isomorphism.
In general, it is a finite cover.)  Let $S_{L}^{*} = (L^{*}\cap
S^{*})^{\circ}$, and let $S_{\Delta}^{*}$ be a subtorus of~$S^{*}$ that is
complementary to~$S_{L}^{*}$ and contains $S^{*}\cap T^{*}$. Let
$S_{\perp}^{*}$ be the image of~$S_{\Delta}^{*}$ under the projection
$L^{*}\times T^{*} \to L^{*}$, and notice that $S_{L}^{*}\cap S_{\perp}^{*} =
e$. Thus, we have
 $$ S^{*} = S_{L}^{*}\times S_{\Delta}^{*} \subset S_{L}^{*} \times
S_{\perp}^{*}\times T^{*} \subset L^{*} \times T^{*}. $$
 Because $\Gamma\langle g\rangle$ is dense in~$G$ and $\Gamma^{\pi} = e$, we
see that $S_{\Delta}^{*}$ finitely covers~$T^{*}$, via the projection
 $S_{\perp}^{*}\times T^{*} \to T^{*}$.

Let $L$ be a semisimple Levi subgroup of~$G$ with $\Ad_{G}L = L^{*}$, and let
$\overline L$ be the corresponding Levi subgroup of~$\overline G$. Since $\Rad
\overline G$ is nilpotent, and $\overline \Gamma\backslash \overline G$ is
faithful, we know that $Z(\overline G)$ is compact (see Cor.~4.6), so some
compact torus $\overline S_{\perp}\subset \overline L$ finitely
covers~$S_{\perp}^{*}$, via the map~$\Ad_{\overline G}$.

By construction, we know that $T^{*}$ is finitely covered by $S_{\Delta}^{*}$,
so the homomorphism $\pi\colon G \to T^{*}$ lifts to a homomorphism $G\to
S_{\Delta}^{*}$. By composing this with the projection $S_{\Delta}^{*} \to
S_{\perp}^{*}$, we obtain a homomorphism $G \to S_{\perp}^{*}$. Then, since
$\overline S_{\perp}$ finitely covers~$S_{\perp}^{*}$, this lifts to a
homomorphism ${\overline \pi}\colon G \to \overline S_{\perp}$.  Note that,
from the definition of~${\overline \pi}$, we have
 $ (\Ad_{\overline G} (x^{{\overline \pi}}) , x^{\pi}) \in S_{\Delta}^{*} $,
 for all $x\in G$. Because $S_{\Delta}^{*}\subset S^{*}$ is contained in the
identity component of the Zariski closure of $\Ad_{G}\langle g\rangle $, which
centralizes~$g$, this implies that 
 \begin{equation} \label{commutes}
 g^{x^{-\pi} x^{-{\overline \pi}}} = g .
 \end{equation}

Replacing $\Gamma$ by a finite-index subgroup, we may assume
$\Gamma^{{\overline \pi}} = e$; hence ${\overline \pi}$ induces a well-defined
map from $\Gamma\backslash G$ to~$\overline S_{\perp}$. Thus, we may define a
homeomorphism
 \begin{equation} \label{psi-def}
 \psi\colon \Gamma\backslash G \to \overline \Gamma\backslash
\overline G \colon x \mapsto x^{\overline \phi}\cdot x^{-{\overline \pi}}. 
 \end{equation}
 Then
 $$ (xg)^{\psi} \={psi-def} (xg)^{\overline \phi}\cdot (xg)^{-{\overline \pi}}
 \={phi-crossed} (x^{\overline \phi} g^{x^{-\pi}\phi})
 (x^{-{\overline \pi}} g^{-{\overline \pi}})
 \={psi-def} x^{\psi} g^{x^{-\pi} x^{-{\overline \pi}} \phi} g^{-{\overline
\pi}}
 \={commutes} x^{\psi} g^{\phi} g^{-{\overline \pi}}
 . $$
  In other words, $\psi$ is affine for~$g$ via~$g^{\phi} g^{-{\overline \pi}}$.
 \end{proof}

\begin{rmk}
 If $G$ is not connected, then, because $G = G^{\circ}\langle g\rangle $, there
is no harm in assuming $G = G^{\circ} \prodsemi \langle g\rangle $, and we may
assume $G^{\circ}$ is simply connected. 

Let $G^{*}$ and $G^{-}$ be the identity components of the Zariski closures of
$\Ad G$ and $\Ad G^{\circ}$, respectively. By replacing $g$ with a
power~$g^{n}$, we may assume $\Ad g \in G^{*}$. Let $S^{+}$ be a maximal
compact torus of the Zariski closure of $\Ad_{G}\langle g\rangle$, and let
$S^{*} = (S^{+}\cap G^{-})^{\circ}$.
 Let $L^{*} T^{+}$ be a reductive Levi subgroup of~$G^{*}$, containing~$S^{+}$,
and define $T^{*} = T^{+}\cap G^{-}$.

 Let $T_{L}T^{+}$ be a maximal compact torus in~$G^{*}$,
containing~$S^{+}$ (where $T_{L}$ is a compact torus in ~$L^{*}$), so
$T_{L}T^{*}$ is a maximal compact torus in~$G^{-}$. Then
$T_{L}T^{*}S^{+} = T_{L}T^{+}$, so there is a subtorus~$Z$ of~$S^{+}$
that is almost complementary to $T_{L}T^{*}$ in $T_{L}T^{+}$.  Assume
for simplicity that $ZT_{L}\cap T^{*} = e$, so there is a natural
projection $G^{*}\to T^{*}$.  Then we may define a homomorphism
$\pi\colon G \to T^{*}$ by composing $\Ad_{G}$ with this projection.

We now construct a semidirect product $G^{*} \prodsemi T^{*}$. Let
 $$ H = \frac{G^{\circ} \prodsemi G^{*}}{\{ (x^{-1}, \Ad x) \mid x
\in [G,G] \}} .$$
 We may assume $\Ad g^n \not\in \Ad G^{\circ}$, for all $n \not= 0$,
for otherwise we could assume $G$~is connected. Therefore, $G \iso
G^{\circ} \langle  \Ad g \rangle $ injects into~$H$. Because 
 $$ [\Ad g, T^{*}] \subset [G^{*},G^{*}] = \Ad [G,G], $$
 we see that $T^{*} \subset N_H(G)$. Therefore, $T^{*}$ is a group of
automorphisms of~$G$, so we may form the semidirect product $G
\prodsemi T^{*}$.

With these definitions of $G^{*}$, $T^{*}$, $S^{*}$, $L^{*}$,
$\pi$, and~$G^{*} \prodsemi T^{*}$ in hand, we may proceed essentially
as in the proof above.
 \end{rmk}

\end{document}